\documentclass[11pt,epsf]{article}
\usepackage{srcltx}
\usepackage{graphicx}
\usepackage{theorem, amsmath, amssymb}
\usepackage{color}
\newtheorem{Theorem}{Theorem}[section]
\newtheorem{Definition}[Theorem]{Definition}
\newtheorem{Proposition}[Theorem]{Proposition}

\newtheorem{Lemma}[Theorem]{Lemma}

\theoremstyle{remark}
\newtheorem{Remark}[Theorem]{Remark}

  {\begin{equation}\left\{\begin{aligned}}%
  {\end{aligned}\right.\end{equation}\ignorespacesafterend}%

\newenvironment{SDE*}%
  {\begin{equation*}\left\{\begin{aligned}}%
  {\end{aligned}\right.\end{equation*}\ignorespacesafterend}%

\begin{document}
\title{
Zero-sum stochastic differential games without the Isaacs condition:\\
random rules of priority and  intermediate Hamiltonians 
}

\author{Daniel Hern\'andez-Hern\'andez \footnote{Research Center for Mathematics (CIMAT), Department of Probability and Statistics, Apartado Postal 402, Guanajuato, GTO, Mexico. E-mail address: dher@cimat.mx. The research of this author was partially supported by CONACYT, under grant 254166.}\; and Mihai S\^{\i}rbu \footnote{University of Texas at Austin,
    Department of Mathematics, 1 University Station C1200, Austin, TX,
    78712.  E-mail address: sirbu@math.utexas.edu. The research of
    this author was supported in part by the National Science
    Foundation under Grant DMS     1517664. Any opinions, findings, and conclusions or recommendations expressed in this material are those of the authors and do not necessarily reflect the views of the National Science Foundation.}}
\maketitle
\begin{abstract} For a zero-sum stochastic game which does not satisfy the  Isaacs condition, we provide a value function representation for an  Isaacs-type  equation whose Hamiltonian lies in between  the lower and upper Hamiltonians, as a convex combination of the two.  For the general case (i.e.  the convex combination is time and  state dependent) our representation amounts to a random change of the rules of the game, to allow each player \emph{at any moment} to see the other player's action or not, \emph{according to a coin toss} with probabilities of heads and tails given by the convex combination appearing in the PDE. If the combination is state independent, then the rules can be set all in advance, in a deterministic way. This means that tossing the coin along the game, or tossing it repeatedly right at the beginning leads to the same value. The representations are asymptotic, over time discretizations. Space discretization is possible as well, leading to similar results.\end{abstract}
\noindent{\bf Keywords:}  stochastic games,  asymptotic  Perron's method, Markov strategies, viscosity solutions, discretization.

\noindent
{\bf Mathematics Subject Classification (2010): }
91A05,  
91A15,  
49L20,  
49L25  

\section{Introduction}
Zero-sum stochastic games of the type
\begin{equation}\label{game}
\sup _{u }\inf _v\mathbb{E}[g(X^{s,x;u,v}_T)],\ \ \ \inf _v\sup _u  \mathbb{E}[g(X^{s,x;u,v}_T)],\ \ \ \ \  \ 
\end{equation}
where the state system is controlled by two opponents $u$ and $v$ as 
\begin{equation}\label{eq:SDE}\left\{
\begin{array}{ll}
 dX_t=b(t,X_t,u_t, v_t)dt+\sigma (t, X_t,u_t, v_t)dW_t,   \\ 
 X_s=x \in \mathbb{R}^d,
 \end{array} \right.
 \end{equation}
have been studied extensively. We recall quickly that 
 the above game is set so that player $v$ pays $g(x)$ to the player $u$ at terminal time $T$, if the state sits at $x$ at time $t=T$.   We mention  just a few representative pieces of work on such zero-sum differential games: the seminal work of Isaacs on deterministic games \cite{isaacs}, some of the early work on stochastic games like \cite{elliott-stochastic-games}, and  the extensive contribution  of the Ekaterinburg school of games (mostly deterministic and  some stochastic)  presented in the monograph  \cite{krasovskii-subbotin-88}. The list is much larger and we do not claim to cover it here. 

It is well known that a (heuristic)  dynamic programming approach to  such  zero sum games leads to the so-called Isaacs equations:
\begin{equation}\label{eq:Isaacs}
\left \{
\begin{array}{ll}
-v_t-H ^{i}(t,x,v_x,v_{xx})=0\ \ \textrm{on}\ [0,T)\times \mathbb{R}^d,\\
v(T,\cdot)=g(\cdot),\ \ \textrm{on}\ \mathbb{R}^d,
\end{array}
\right.
\end{equation}
where
\begin{equation}\label{lower-h}
H^-(t,x,p,M)\triangleq \sup_{u \in U} \inf _{v\in V}\left[b(t,x,u,v)\cdot p+\frac{1}{2}Tr \left (\sigma(t,x,u,v)\sigma(t,x,u,v)^T M\right)\right],\  0\leq t\leq T, \; x\in \mathbb{R}^d,
\end{equation}

\begin{equation}\label{upper-h}
H^+(t,x,p,M)\triangleq \inf _{v\in V}\sup_{u \in U}\left[b(t,x,u,v)\cdot p+\frac{1}{2}Tr\left (\sigma(t,x,u,v)\sigma(t,x,u,v)^T M\right)\right],\  0\leq t\leq T, \; x\in \mathbb{R}^d.
\end{equation}

For $i=-$ we call it  the lower Isaacs equation, and for $i=+$ we call it the upper Isaacs equation. 
The precise definition of the games varies in the literature (as there is no canonical model), from the Elliott-Kalton formulation in \cite{fs}  to the feedback formulation in \cite{krasovskii-subbotin-88} which is  used in similar ways in \cite{fleming-hernandez2}, \cite{fleming-hernandez2-2}, \cite{sirbu} or \cite{sirbu-4}. However, in any  of these formulations, the intuition behind  the (solution of the) lower equation is that player $v$ has informational priority over $u$ . By \emph{informational priority}  we mean that  $v$ can see, \emph{in real time}, the  action of $u$ before  choosing his/her own action $v$. The upper Isaacs equation has, obviously, the exact opposite meaning. 
If the 
 so-called Isaacs condition holds, i.e.
$$H^-=H^+ \triangleq H,$$
it does not matter who has informational priority. Even if both players can (symmetrically) only observe the state and nothing else, the game still have a value (see \cite{sirbu}).

The main goal of the present paper is to provide a meaningful  interpretation/representation for the (unique viscosity) solution  of the 
 equation 
\begin{equation}\label{eq:Isaacs-p}
\left \{
\begin{array}{ll}
-v_t-H ^p (t,x,v_x,v_{xx})=0\ \ \textrm{on}\ [0,T)\times \mathbb{R}^d,\\
v(T,\cdot)=g(\cdot),\ \ \textrm{on}\ \mathbb{R}^d,
\end{array}
\right.
\end{equation}
where
$$H^p(t,x,p,M)=p(t,x) H^-(t,x,p,M)+  (1-p(t,x))H^+(t,x,p,M),$$
for some fixed $p:[0,T]\times \mathbb{R}^d\rightarrow [0,1]$.
Such representation should naturally come as a value obtained from playing the formal game \eqref{game}-\eqref{eq:SDE}. This means that we are looking for a model of rules of the game and strategies  such that the resulting game has  a  value and the value is equal to the solution of \eqref{eq:Isaacs-p}. In order to model the strategies of the players  we follow here the line of work  \cite{krasovskii-subbotin-88}, \cite{fleming-hernandez2}, \cite{fleming-hernandez2-2}, \cite{sirbu} or \cite{sirbu-4} that  uses feedback strategies, rather than the idea of Elliott-Kalton strategies in \cite{ek}, \cite{fs}.

\noindent{\bf The Rules of the Game:} in words, the (heuristic) rules of a game with  value $v^p$ are\begin{enumerate}
\item {\bf feedback strategies:} at each time $t$, both players observe the full past of the state $X|_{[0,t]}$ ,
\item {\bf randomization:} at each such time $t$, if the state sits at $X_t=x$ a   coin  is tossed, with probabilities $\mathbb{P}(H)=p(t,x)$, $\mathbb{P}(T)=1-p(t,x)$,
\item {\bf priority rule:} if the coin turns $H$, then at that particular time, $v$ sees $u$  first,in real time.  If the coin turns  $T$,   then $u$ sees $v$.
\end{enumerate}

\begin{Remark}\label{rem:cont-game}Formally, if the  Isaacs equation \eqref{eq:Isaacs-p} has a smooth solution $v^p$, playing a continuous-time game with value $v^p$ is  nothing else then a "local game" over $(u,v)$ (that  should result in the value $v^p$).  
\end{Remark}
Remark \ref{rem:cont-game}, as well as the heuristic rules of the game described above,  are all justified by  the very simple one-period observation in Subsection \ref{intro:sub} below. The rest of the paper is dedicated to putting the game with such rules on sound theoretical foundation, and proving the (asymptotic) value property. When the coin toss does not depend on the state, we actually have a deterministic approximation, in the spirit of \cite{kaise-sheu}. 
The  heuristic connection  between the two cases, the  random rule of priority or  the (quick) deterministic change for  the rule of priority,  fits in the well known narrative where, in a  continuous-time problem, one can basically simulate randomization from the realization of a single path.
  Since the  technical proofs  for the deterministic case are a little bit more involved, we present this result (and its proof) first.
  
Technically, the proofs are mainly based on a modification of Perron's method introduced in \cite{sirbu-4}. It should be noted that while  Perron scheme  is related to the  monotone scheme of \cite{barles-souganidis} (related in turn  to the method of relaxed semi-limits of   \cite{barles-perthame-1}), the version  of Perron used here is based on \emph{probabilistic} arguments.

To the best of our knowledge, the only paper studying a similar representation problem so far is the work \cite{kaise-sheu}, for deterministic games. For the  case when the problem is \emph{ autonomous} with respect to time \emph{and}  $p$ is a \emph{constant} (probability)  they use a Trotter-Kato approximation idea to show the following: if one denotes by $S^-(t)$ and $S^+(t)$ the non-linear (Nisio-type, see \cite{MR945913}) semigroups  associated with the two games, i.e.
$v^-(t,\cdot) =S^-(T-t)g$ and $v^+(t,\cdot) =S^+(T-t)g$ are (the viscosity) solutions of the two Isaacs equations \eqref{eq:Isaacs} for $i=-$ and $i=+$  given by \cite{fs}, then we have the (analytic) result 
\begin{equation}\label{eq:ks}v^p(t,\cdot)=\lim _{n\rightarrow \infty} \Big (S^-(p (T-t)/n)S^+(p (T-t)/n)\Big)^ng.
\end{equation} 
The right-hand side above  \emph{might} be interpreted as a time-discretized game,  but the precise definition of  a game with such a value is actually  not necessarily possible. Therefore,  the result is mostly an analytical asymptotic representation. In addition, the semigroups $S^-$ and $S^+$ are constructed in \cite{fs} using non-symmetric games which play Elliott-Kalton strategies vs. open-loop controls, introduced for the deterministic framework in \cite{ek}. 
As  observed in \cite{sirbu-4}[Remark 2.1], these definitions do not  always have the  genuine meaning of some  values of non symmetric games. This is why it may not be possible to  rewrite the right-hand-side in \eqref{eq:ks} as an inf/sup or sup/inf.

While our work below in Subsection \ref{deterministic} is, obviously, related to \cite{kaise-sheu}, it does not really overlap, either conceptually or technically.  To begin with, for us  $v^p$ is a limit of a genuine (sup/inf  or inf/sup) value of a game, discretized over time (compare to the comments above). Our state system is time dependent, and, more importantly, we allow for the probability $p$ to be non-constant.
The time discretization we use is different, as we do not allow for actions to be changed in between the times on the grid, and the time grid does not have to be as in \cite{kaise-sheu} but only satisfy what we call an asymptotic density property.  This relates to the possibility of randomization of rules (as in the general Subsection \ref{random}) but where randomization of rules is performed in advance.

Most importantly, our most general result in Subsection \ref{random} provides a completely different conceptual point of view, as a game with \emph{random rules}, were \emph{a time and state dependent coin} is tossed to decide which player has priority. Such a result, even at a formal level, does not exist in the literature, to the best of our knowledge,  even in the deterministic case.

Technically, our  work (both parts) is   rather different from \cite{kaise-sheu}, being based on a Perron scheme introduced in \cite{sirbu-4}. One of the benefits is that we do \emph{not} need to assume any prior results from the theory of zero-sum games (the results in \cite{fs} are used in \cite{kaise-sheu}), but we can build the  whole model and perform the analysis in a  self-contained way. 
\subsection{One-period games with random rules}\label{intro:sub}
In the general case of Subsection \ref{random}, our dynamic game will be played using \emph{locally} a simple idea that is best explained in one-period.   This  is easy enough so that  genuine proofs are unnecessary. The notation  $u,v$ for the actions  of the players overlaps with the dynamic case, but this is on purpose, to  facilitate the presentation.

Assume we have a (bounded) function  $f:U\times V\rightarrow \mathbb{R}$ and we  define
$$f^-\triangleq \sup_{u \in U} \inf _{v\in V}f(u,v)\leq  \inf _{v\in V} \sup_{u \in U}f(u,v)\triangleq f^+.$$
We can assign to $f^-$ the meaning of a true value/saddle point for the game where $v$ sees $u$ and similarly for $f^+$. More precisely, if we denote by $\beta$ the possible responses for $v$ that sees $u$, i.e. functions $\beta:U\rightarrow V$  then 
$$f^-=\sup_{u} \inf _{\beta }f(u,\beta (u))=  \inf _{\beta } \sup_{u }f(u,\beta (u)).$$
Similarly, if  $u$ sees the action of $v$,  considering all maps  $\alpha :V\rightarrow U$ then 
$$f^+=\sup_{\alpha } \inf _{v\in V}f(\alpha (v),v)= \inf _{v\in V} \sup_{\alpha }f(\alpha (v),v).$$
In the spirit of  the work on dynamic games by Krasovskii-Subbotin  \cite{krasovskii-subbotin-88} we  call  $\alpha$ and $\beta$ above  "counter-strategies" in this static game.
Obviously, the saddle points above (if they exist) are found by  solving a sequential optimization problem, i.e. finding  (for $f^-$, for example) first 
$$\beta ^* (u)=\arg \min_{v} f(u,v), \ \textrm{then}\ \  u^*=\arg\max _u \inf_{v} f(u,v)=\arg\max_{u} f(u, \beta ^* (u)).$$
We similarly find $\alpha ^*$ and $v^*$ for $f^+$  (if they exist, if not, we can find them approximately).

Now, the interpretation of $$f^p\triangleq p f^-+(1-p) f^+$$ for some $0\leq p\leq 1$ is that of the value of a game played by the following rules
\begin{enumerate}
\item a coin is tossed, with probability of $H$ equal to $p$
\item if $H$ shows up then, by the rule of the game, player $v$ sees the control of player $u$
\item if $T$ shows up, then, by the rule of the game, player $u$ sees the control of player $v$.
\end{enumerate}
This can also be thought of as: \emph{prior} to the coin toss, player $u$ chooses a pair $(u,\alpha)$ and player $v$ chooses a pair $(v,\beta)$.
Then \begin{enumerate}
\item if $H$ shows up then $u$ is played against $\beta$
\item if $T$ shows up, then $\alpha$ is played against $v$.
\end{enumerate}
In other words, we have the quite obvious representation
$$f^p=\sup_{(u, \alpha)} \inf _{(v,\beta)}  \big [p f(u,\beta (u))+(1-p)f(\alpha (v),v) \big] = 
\inf _{(v,\beta)}\sup_{(u, \alpha)}   \big [p f(u,\beta (u))+(1-p)f(\alpha (v),v) \big].
$$
In addition, the pairs $(u^*, \alpha ^*)$ and $(v^*, \beta ^*)$  found above (again, if such  pairs exist, if not they exist in an approximate sense) represent a saddle point for the game with value $f^p$.
Once again, the general dynamic game with a value $v^p$ in Subsection \ref{random} will be played using this idea continuously at any time $t$.
\section{Set-up and main results}  We now turn to the precise conditions on the state equation and pay-off, and leave the precise definition of the  probability space, as well as strategies for slightly later.
We have here a  differential game with two players.
The first player's actions belong to the  compact metric space ($U, d_u)$
(usually $ U\subset \mathbb{R}^k$). The second player's actions belong also to a  compact metric space $(V,d_v)$. We assume that the state belongs to  $\mathbb{R}^d$ (an   open domain $\mathcal{O}\subset \mathbb{R}^d$ as in \cite{bs-3} is also feasible).
The coefficients  $b:[0,T] \times \mathbb{R}^d \times U \times V \to \mathbb{R}^d$ and $\sigma:[0,T]\times\mathbb{R}^d\times U \times V \to \mathbb{M}_{d,d'}$ of the state equations  satisfy some usual assumptions (see \cite{sirbu} or \cite{sirbu-4}).

\pagebreak

\noindent {\bf Assumptions on the  state system:}   $b$ and $\sigma$ are 

\begin{enumerate}
\item continuous on $[0,T] \times \mathbb{R}^d \times U \times V$,
\item uniformly locally Lipschitz in the state variable, i.e. for each $K<\infty$ there exists $L(K)<\infty $ such that
$$|b(t,x,u,v)-b(t,y,u,v)|+|\sigma (t,x,u,v)-\sigma (t,y,u,v)|\leq L(K)|x-y|, \ \ \forall |x|,|y|\leq K, \ \ \forall t,u,v,$$
\item have linear growth, i.e. there exists $C<\infty$ such that 
$$|b(t,x,u,v)|+|\sigma (t,x,u,v)|\leq C(1+|x|),\ \ \forall t,x,u,v.$$

\end{enumerate}
In addition, we have the following

\noindent {\bf Assumption on the pay-off:} the function $g:\mathbb{R}^d\rightarrow \mathbb{R}$ is continuous and bounded.

 The state of the system is governed by equation \eqref{eq:SDE}, if the game  starts at an initial time $s$ at some  position $x$. In order to obtain the desired representation for the solution of the Isaacs-type equation \eqref{eq:Isaacs-p}, we will allow, as mentioned in the introduction, for changes of rules of the game, as time evolves. We do so by discretizing time.  State discretization is also possible. 
 For a fixed $0\leq s\leq T$, we will usually denote by $\Delta$ a time partition
$$s=t_0\leq t_1\leq ...\leq t_n=T, \ \ {for \ some} \ n.$$
and by $\mathcal{D}(s)$ the collection of all partitions of the interval $[s,T]$. 

\subsection{Deterministic rules of priority: state-independent intermediate Hamiltonians}\label{deterministic}  Fix $(s,x)$. Assume here that the process $W$ is a $d'$-dimensional Brownian motion on a \emph{fixed} probability space
 $(\Omega, \mathcal{F}, (\mathcal{F}_t)_{s\leq t\leq T},\mathbb{P})$ (which can, actually, depend on $(s,x)$). The filtration $(\mathcal{F}_t)_{s\leq t\leq T}$ satisfies the usual conditions and may be larger than the natural filtration
$\mathcal{F}^W_t\triangleq \sigma \{W_u;\ s\leq u\leq t\}\vee \mathcal{N}(\mathbb{P}, \mathcal{F}).$

 Consider  now  a time partition $\Delta=(t_0, t_1, \dots, t_n) \in \mathcal{D}(s)$ of $[s,T]$. Together with this partition we have a sequence   $\xi=\xi_1, \dots, \xi_n$ of "marks" valued as  $0$ or $1$.  We denote by $\mathcal{M}(\Delta)$ the collection of all possible sequence of marks (note that we have as many marks as time intervals in the partition). 
 
 The rules of the game are set, deterministically, in advance, and the players cannot change actions in between the discrete times on $\Delta$. If the mark $\xi_k$ is equal to zero, then, on the interval $[t_{k-1}, t_k]$ player u first sees player v's action, i.e player u chooses (over this time interval) some $\alpha$
"counterstrategy " depending on 
\begin{enumerate}
\item the whole past of the state up to $t_{k-1}$
\item current value of the other players' action 
\end{enumerate} 
The other player only choose a strategy $b$ depending on the past of the state, up to $t_{k-1}$. In other words, if $\xi_k=0$ the first player $u$ uses over  $[t_{k-1}, t_k]$ a
$\alpha _k(X|_{[0, t_{k-1}]}, v)$ and player $v$ uses a
$b _k(X|_{[0, t_{k-1}]})$.
If $\xi_k=1$ then the roles are reversed, and player $u$ uses a strategy $a _k(X|_{0, t_{k-1}]})$ while player $v$ uses a counter-strategy  $\beta _k(X|_{[0, t_{k-1}]}, v)$. These can  all  me bade "Markov" (with some special care)  or general dependent on the whole past, or may even  depend, \emph{explicitly}, on the actions  of the opposing player over the intervals $[t_0, t_1],\dots [t_{k-1}, t_{k-1}]$. 
The marks are deterministic, i.e. are known in advance when the game is started.
With these rules, we can define a value of the game, depending on the partition and the marks
$\xi =(\xi_1, \dots, \xi_n).$
Before that, for a  fixed $s$, we denote by $C[s,T]$ the set of continuous paths $C([s,T]:\mathbb{R}^d)$. A generic path will be denoted by $y$ or $y(\cdot)$.  We also denote by
$\mathcal{B}_t$ the  raw filtration generated by the paths up to time $t$, i.e.
$\mathcal{B}_t=\sigma \{y(u): s\leq u\leq t\}.$ 

\begin{Definition}[Discrete strategies, deterministic change of rules]  Fix $s$, together with  the time partition $\Delta \in \mathcal{D}(s)$ and the deterministic change of rules  $\xi \in \mathcal{M}(\Delta)$. 
\begin{enumerate} \item A  full strategy for the player $u$ is a sequence
$$ \mathbf{u}= (u_k), k=1,2,...n,$$ where 
\begin{enumerate}
\item if $\xi_k=1$ then 
$u_k=a_k: C[s,T]\times V^{k-1} \rightarrow U $ and
\item if $\xi_k=0$ then $  u_k=\alpha _k: C[s,T]\times V^k \rightarrow U$.
\end{enumerate}
Here, $a_k$ is  
$\mathcal{B}_{t_{k-1}}\otimes\mathcal{B}(V^{k-1})/\mathcal{B}(U)$ measurable and  $\alpha _k$ is  $\mathcal{B}_{t_{k-1}}\otimes\mathcal{B}(V^k)/\mathcal{B}(U)$  measurable.   We denote by
$\mathcal{A}^{s,\Delta, \xi}$ the set of these strategies. Player $u$  chooses the action $u_k$ at time $t_{k-1}$  and holds  it until $t_k$.  This depends on the whole past of the state up to $t_{k-1}$, the actions $v_1,\dots, v_{k-1}$ of the opponent on the intervals $[t_0,t_1], \dots , [t_{k-2}, t_{k-1}]$ and, depending on the mark $\xi_k$ on the value $v_k$ or not.

\item  A  feedback strategy for the first player is a sequence
$$ \mathbf{u}= (u_k), k=1,2,...n,$$ where

\begin{enumerate}
\item if $\xi_k=1$ then 
$u_k=a_k: C[s,T] \rightarrow U $ and
\item if $\xi_k=0$ then $  u_k=\alpha_k: C[s,T]\times V \rightarrow U$.
\end{enumerate}
where $a_k$ is $\mathcal{B}_{t_{k-1}}/\mathcal{B}(U)$ and $\alpha _k$ is  $\mathcal{B}_{t_{k-1}}\otimes\mathcal{B}(V)/\mathcal{B}(U)$  measurable.
The meaning is again clear.  We denote by $\mathcal{A}_F^{s,\Delta,\xi}$ the set of these strategies
\item  Markov strategy for the first player: the definition here is a little bit different, because we allow for the player to change actions fewer times than the marks change. 
More precisely,  for  the fixed $\Delta: t_0<\dots t_n=T$  we assume that there is a sub-grid
$$s=r_0 <\dots <r_{I(n)}, \ \ \ \ \ r_i= t_{l(i)}, \ i=0, \dots,  I(n)\leq n,$$
such that, at any of these times $r_{i-1}=t_{l(i-1)}$, the player chooses a pair $(a_i, \alpha_i)$ such that 
$a_i$ is $\mathcal{B}(\mathbb{R}^d) /\mathcal{B}(U)$ and $ \alpha_i$ is  $\mathcal{B}(\mathbb{R}^d)\otimes\mathcal{B}(V)/\mathcal{B}(U)$  measurable and holds the action decided based on where  the state is at time $r_{i-1}=t_{l(i-1)}$ until the next time on the sub-grid, i.e. $r_i=t_l(i)$.  The marks may change  in between $r_{i-1}$ and $r_i$ but the action of player $u$ only changes if the mark is $\xi_k=0$ and the action of the opponent changes.
In other words, fixed the double  sequence $(a_i,\alpha _i)_{i=1}^{I(n)}$, the actions of player $u$ can be represented as a strategy in item 1 (with some abuse of notation) by
$$ \mathbf{u}= (u_k), k=1,2,...n,$$ where, 
if $l(i-1)<k\leq l(i)$ then 
$$u_k(X_{[0, t_{k-1}]})=\left \{\begin{array}{ll}
a_i (X_{t_{l(i-1)}}), & \xi _k=1\\
\alpha _i (X_{t_{l(i-1)}}, \cdot )\,  & \xi _k=0. \end{array}\right. $$
We denote by
$\mathcal{A}_M ^{s,\Delta,\xi}$ the set of these strategies. \end{enumerate}
\end{Definition}
We have
$$\mathcal{A}_M^{s,\Delta,\xi}\subset  \mathcal{A}_F^{s,\Delta,\xi}\subset \mathcal{A}^{s,\Delta,\xi}.$$
The first inclusion may not seem obvious, but, at a second glance, it is clear that dependence on older past is OK. 
We define in a similar  manner (but symmetrically opposite, depending on the value of the marks) the (sets of) strategies for the player $v$, 
$\mathcal{B}_M^{s,\Delta,\xi}\subset  \mathcal{B}_F^{s,\Delta,\xi}\subset \mathcal{B}^{s,\Delta,\xi}.$
In the definition of strategies $\mathbf{v}$ the role of the value of the marks $\xi$ is reversed, in an obvious way. 
\begin{Proposition} Fix $s, x$ as well as $\Delta \in \mathcal{D}(s)$, $\xi \in\mathcal{M}(\Delta)$.  Choose  $\mathbf{u}\in  \mathcal{A}^{s,\Delta, \xi }$,
$\mathbf{v}\in \mathcal{B}^{s,\Delta, \xi } $.  There exists a unique strong solution of the state system \eqref{eq:SDE}, denoted by 
$(X^{s,x;\mathbf{u}, \mathbf{v}}_t)_{s\leq t\leq T}$, such that
$$X_t\in \mathcal{F}^W_t \subset \mathcal{F}_t , \  s\leq t\leq T.$$

\end{Proposition}
Proof: the proof is done step-by-step over intervals $[t_{k-1}, t_k]$ as in \cite{krasovskii-subbotin-88} or \cite{sirbu}.  Actually, since the initial condition is deterministic, the (uniform) linear growth condition insures that the solution is square integrable. $\diamond$

 With  the notation
$J(s,x, \Delta, \xi;\mathbf{u}, \mathbf{v})\triangleq \mathbb{E}\left [g \left ( X^{s,x;\mathbf{u}, \mathbf{v}}_T \right) \right]$,  we can now define the value functions:
\begin{equation}
\begin{split}
V_M^-(s,x;\Delta,\xi) & \triangleq \sup _{\mathbf{u}\in \mathcal{A}_M^{s,\Delta,\xi }} \inf _{\mathbf{v}\in \mathcal{B}^{s,\Delta ,\xi}} J(s,x, \Delta, \xi;\mathbf{u}, \mathbf{v}) \leq  \\
V^-(s,x;\Delta,\xi ) & \triangleq \sup _{\mathbf{u}\in \mathcal{A}^{s,\Delta,\xi}}  \inf _{\mathbf{v}\in \mathcal{B}^{s,\Delta,\xi}}  J(s,x, \Delta, \xi;\mathbf{u}, \mathbf{v}) \leq \\
 V^+(s,x;\Delta, \xi ) & \triangleq \inf _{\mathbf{v}\in \mathcal{B}^{s,\Delta, \xi}} \sup _{\mathbf{u}\in \mathcal{A}^{s,\Delta ,\xi}} J(s,x,\Delta,\xi;\mathbf{u}, \mathbf{v}) \leq \\
V_M^+(s,x;\Delta,\xi) & \triangleq \inf _{\mathbf{v}\in \mathcal{B}_M^{s,\Delta,\xi}} \sup _{\mathbf{u}\in \mathcal{A}^{s,\Delta,\xi}} J(s,x,\Delta, \xi;\mathbf{u}, \mathbf{v}).\end{split} \end{equation}
In the chain of inequalities above we could have defined an additional intermediate layer of value functions, playing Markov strategies vs. feedback strategies.
Next theorem  is a (deterministic) asymptotic result, somewhat similar to \cite{kaise-sheu}, but different: time intervals are not ordered, and the actions cannot be continuously changed in between the times on the grid. Also, we have a stochastic game for which we find asymptotic values of inf/sup and sup/inf type.
\begin{Theorem}\label{main-deterministic} Assume the function 
$p:[0,T]\rightarrow [0,1]$ is continuous. For each initial time $s$
consider  a sequence of partitions  $\Delta ^{n,s}\in \mathcal{D}(s)$, together with marks  $\xi ^{n,s} \in \mathcal{M}(\Delta ^{n,s}))$ such that 
$(\Delta ^{n,s}, \xi^{n,s})$ has {\bf uniform} asymptotic density $p$. By this, we mean that for 
 any  $\varepsilon>0$ there exists $n_0(\varepsilon)$ (in this formulation we assume that $n_0$  {\rm does not} depend on $s$) large enough, such that for any  $n\geq n_0(\varepsilon)$ we have that, inside the partition $$\Delta ^{n,s}=s=t_0<t_1....<t_{N(n)}=T,$$  we can find a sub-division
$$s=r_0<r_1<...r_{I(n)}=T,\ \  \ \ r_i=t_{l(i)},\ \ l(0)=0 <l(1)<...<l(I(n))=N(n), $$
such that
$$r_i-r_{i-1}\leq \varepsilon,  \forall \ i=1,...,I(n), $$
and
$$\left |\frac { \sum _{k=l(i-1)+1} ^{l(i)}(t_k -t_{k-1}) \xi_k }{t_{l(i)}-t_{l(i-1)}}-p(t_{l(i-1)})\right | \leq \varepsilon,\ \ \ 
 \forall \ i=1,...,I(n)\leq N(n).$$
Then the Isaacs-type equation \eqref{eq:Isaacs-p} has a unique bounded continuous viscosity  solution $v$,   which is  the limit  (uniform on compacts sets)  of  both $V^-_M(\cdot, \cdot;\Delta ^{n,\cdot}, \xi ^{n, \cdot} )$ and $V^+_M (\cdot, \cdot;\Delta ^{n, \cdot}, \xi^{n, \cdot}) $.
More precisely, for any $K< \infty$ and $\varepsilon >0$, there exists  $n_0=n(K, \varepsilon)$ such that
$$\forall s, \forall |x|\leq K,\ \ \forall n\geq n_0,$$ we have
\begin{equation}\label{main-det-precise}V_M^+(s,x;\Delta ^{n,s}, \xi ^{n,s})-\varepsilon \leq v(s,x)\leq V_M^-(s,x;\Delta ^{n,s},\xi^{n,s})+\varepsilon.
\end{equation}
\end{Theorem}

\subsection{Random rules of priority: state-dependent intermediate Hamiltonians}\label{random}
We study here the general case, where the players play, locally, a game with random rules described in Subsection \ref{intro:sub}. Assume that the convex combination parameter is a continuous function
$$p:[0,T]\times \mathbb{R}^d\rightarrow [0,1].$$ Fix the initial time $s$.  The process $W$ is a $d'$-dimensional Brownian motion on the  \emph{fixed}  filtered probability space 
 $(\Omega, \mathcal{F},( \mathcal{F}_t)_{s\leq t\leq T} \mathbb{P})$.
 The filtration satisfies the usual conditions.
In order to perform independent coin tosses, we  assume that the space $\Omega$ accommodates  the (infinite) sequence of independent standard normals
$$\eta =(\eta _k)_{k=1,2,\dots}.$$
They are independent of the Brownian motion, and, since $W$ comes with a possibly larger filtration, we actually assume that, 
under $\mathbb{P}$,  the whole sequence  $\eta$ is independent of $\mathcal{F}_T$ (and is i.i.d normal, as mentioned).
We denote by $\Phi$ the c.d.f. of a standard normal, so that, for each $\eta_k$ above we have
$$\mathbb{P}(\Phi ^{-1}(\eta _k)\leq \lambda)=\lambda,\ \ 0\leq \lambda \leq 1.$$
Fix now the initial position $x$. Fix, in addition, as above, a partition $\Delta$ of the time interval. 
The game is discretised in time over $\Delta$ and the rules are decided by the realisations of the first $\eta_i$'s (as many as time intervals). 
Instead of considering an infinite sequence of i.i.d. normals, one could even assume that the probability space $(\Omega, \mathcal{F},( \mathcal{F}_t)_{s\leq t\leq T}, \mathbb{P})$, the Brownian motion $W$ and the first $n$ coin tosses $\eta _1, ...\eta_n$ actually depend on $(s,x)$ and $\Delta$.  We could also consider directly  a sequence of i.i.d. uniform on $(0,1)$ instead of $ (\Phi ^{-1}(\eta _k))_{k=1,2,\dots}$.

\noindent {\bf Main rules of the game:}
\begin{enumerate}
\item {\bf Random priority:} at each time $t_{k-1}$ on the grid, depending on where the state $X_{t_{k-1}}$ is located, a coin is tossed with probability $p(t_{k-1}, X_{t_{k-1}})$ and then players $u$, $v$  choose their actions with a priority rule described in the Subsection \ref{intro:sub}  of the Introduction. The coin toss is simulated using the binary random variable
$$1_{\{\Phi ^{-1}(\eta _{k})\leq p(t_{k-1}, X_{t_{k-1}}\}},$$
so the value of $\eta _{k}$ decides which player has priority over the interval $[t_{k-1}, t_k].$
\item {\bf No action outside the grid:} the controls chosen at $t_{k-1}$ have to be kept  \emph{constant} until the next time on the grid, $t_{k}$.
\end{enumerate}
This means that, at time $t_{k-1}$ player $u$ chooses a couple $(u_{k}, \alpha _{k})$ and player $v$ chooses $(v_{k}, \beta _{k})$,   as in  Subsection \ref{intro:sub}  of the Introduction.  
\begin{enumerate}
\item If $1_{\{\Phi ^{-1}(\eta _{k})\leq p(t_{k-1}, X_{t_{k-1}})\}}=1$ then $u_{k}$ is played against $\beta _{k}$  until $t_k$ and, \item if 
$1_{\{\Phi ^{-1}(\eta _{k}\leq p(t_{k-1}, X_{t_{k-1}})\}}=0$ then $\alpha _{k}$ is played against $v _{k}$ at time $t_{k-1}$ and then kept constant until $t_{k}$.
\end{enumerate}
The set-up of the game has to also  take into account the additional (to the coin toss) information  that the two players have at time $t_{k-1}$, in order to choose the two couples $(u_k, \alpha _k)$,  $(v_k, \beta _k)$. Therefore, we need the

\noindent {\bf Additional informational structure:}  there are three symmetric possibilities where (both) players 
\begin{enumerate}
\item at time $t_{k-1}$ can see the whole past of the  state up to that time \emph{and} all controls used by the opposing player in the past,
\item at time $t_{k-1}$  can see the past of the state up to time $t_{k-1}$, but not the past controls,
\item at time $t_{k-1}$ players can see only the state at time $t_{k-1}$.

\end{enumerate}
We use the same notation for  $C[s,T]$ as  the set of continuous paths $C([s,T]:\mathbb{R}^d)$,   generic paths are y $y$ or $y(\cdot)$ and 
$\mathcal{B}_t$ is  the  raw filtration generated by the paths up to time $t$, i.e.
$\mathcal{B}_t=\sigma \{y(u): s\leq u\leq t\}$.

\begin{Definition}[Discrete strategies]  Fix $s$ and $\Delta \in \mathcal{D}(s)$ 
where $\Delta$ is given as 
$$ s=t_0<t_1<\dots<t_n=T\ \ \textrm{for\ some\ } n.$$ 
\begin{enumerate} \item A  full strategy for the first player is a sequence
$$ \mathbf{u}= (a_k, \alpha _k), k=1,2,...n,$$ where
$$a_k: C[s,T]\times V^{k-1} \rightarrow U \ \textrm{and}\ \ \alpha_k: C[s,T]\times V^k \rightarrow U$$
are 
$ \mathcal{B}_{t_{k-1}}\otimes\mathcal{B}(V^{k-1})/\mathcal{B}(U)$ and, respectively, $ \mathcal{B}_{t_{k-1}}\otimes\mathcal{B}(V^k)/\mathcal{B}(U)$  measurable.

In words, the  $k$-th decision of player $u$ made  at time $t_{k-1}$ depends on the  history of the state process up to $t_{k-1}$, and history of the actions of $v$ over intervals $[t_0, t_1], ....[t_{k-2}, t_{k-1}]$, as well as the $k$ coin toss.  We denote by
$\mathcal{A}^{s,\Delta}$ the set of these strategies.

\item  A  feedback strategy for the first player is a sequence
$$ \mathbf{u}= (a_k, \alpha _k), k=1,2,...n,$$ where
$$a_k: C[s,T]\rightarrow U \ \textrm{and}\ \ \alpha_k: C[s,T]\times V\rightarrow U$$
are 
$ \mathcal{B}_{t_{k-1}}/\mathcal{B}(U)$ and, respectively, $ \mathcal{B}_{t_{k-1}}\otimes\mathcal{B}(V)/\mathcal{B}(U)$  measurable.
In words, this is a strategy as above but  the history of the actions of $v$ over intervals $[t_0, t_1], ....[t_{k-2}, t_{k-1}]$ is not taken into account (explicitly) for the $k$'th decision of player $u$ at time $t_{k-1}$. We denote by $\mathcal{A}_F^{s,\Delta}$ the set of these strategies
\item  A  Markov strategy for the first player is a sequence
$$ \mathbf{u}= (a_k, \alpha _k), k=1,2,...n,$$ where
$$a_k: \mathbb{R}^d \rightarrow U \ \textrm{and}\ \ \alpha_k: \mathbb{R}^d \times V\rightarrow U$$
are 
$ \mathcal{B}(\mathbb{R}^d) /\mathcal{B}(U)$ and, respectively, $ \mathcal{B}(\mathbb{R}^d)\otimes\mathcal{B}(V)/\mathcal{B}(U)$  measurable.

In this case, the decision $k$'th  of player $u$, made at time $t_{k-1}$, depends only on the position $X_{t_{k-1}}\in \mathbb{R}^d$, and, depending on  how the $k'$th toss turns up, the other opponent's action at time $t_{k-1}$ or not.
We denote by
$\mathcal{A}_M ^{s,\Delta}$ the set of these strategies. \end{enumerate}
\end{Definition}
Obviously, 
$$\mathcal{A}_M^{s,\Delta}\subset  \mathcal{A}_F^{s,\Delta}\subset \mathcal{A}^{s,\Delta}.$$
We define in an identical manner similar (sets of) strategies for the player $v$, 
$$\mathcal{B}_M^{s,\Delta}\subset  \mathcal{B}_F^{s,\Delta}\subset \mathcal{B}^{s,\Delta}.$$
\begin{Remark} The time-discretization, together with the  assumption that actions of the players are kept constant in between the times on the grid $\Delta$ allow for a simple definition of full strategies. In our previous work \cite{sirbu} and \cite{sirbu-4}, we only used the idea of feedback strategies (or counterstrategies). However,  it is well know from discrete-time game theory, that the idea of full strategy should not lead to different value. We consider it here for completeness. From the intuition point of view, it appears also rather clear that observing the full past of the state allows one to recover all the information needed about the past controls of the opponent. The really important piece of information is the control of the opponent \emph{in real time} at the time of a possible change (on the grid). Modelling precisely how this  is decided randomly is the goal of the paper.
\end{Remark}
\begin{Proposition} Fix $s, x, \Delta$. Fix $\mathbf{u}\in  \mathcal{A}^{s,\Delta}$,
$\mathbf{v}\in \mathcal{B}^{s,\Delta} $.  There exists a unique strong solution of the state system \eqref{eq:SDE}, denoted by 
$(X^{s,x;\mathbf{u}, \mathbf{v}}_t)_{s\leq t\leq T}$, such that
$$X_t\in \mathcal{F}^W_t  \vee \sigma (\eta _1, \dots , \eta _k)\subset \mathcal{F}_t  \vee \sigma (\eta _1, \dots , \eta _k), \ \ t_{k-1}\leq t\leq t_k.$$

\end{Proposition}
Proof: the proof is based, again, on  the  "step-by-step"  scheme as in \cite{krasovskii-subbotin-88}. In other words, one solves the system, iteratively over intervals $[t_{k-1}, t_k]$, noticing that the realization of the strategies  $(a_k, \alpha _k)$ and  $(b_k, \beta _k)$ at time $t_{k-1}$ are kept constant over $[t_{k-1}, t_k]$. The only additional thing needed is the coin toss $\eta_k$, which is independent of $W$.
$\diamond$

With  the notation
$J(s,x;\mathbf{u}, \mathbf{v})\triangleq \mathbb{E}\left [g \left ( X^{s,x;\mathbf{u}, \mathbf{v}}_T \right) \right],$
we  define the value functions:
\begin{equation}
\begin{split}
V_M^-(s,x;\Delta) & \triangleq \sup _{\mathbf{u}\in \mathcal{A}_M^{s,\Delta}} \inf_{\mathbf{v}\in \mathcal{B}^{s,\Delta}} J(s,x;\mathbf{u}, \mathbf{v}) \leq 
V^-(s,x;\Delta)\triangleq \sup _{\mathbf{u}\in \mathcal{A}^{s,\Delta}}  \inf _{\mathbf{v}\in \mathcal{B}^{s,\Delta}}  J(s,x;\mathbf{u}, \mathbf{v}) \leq \\
\leq V^+(s,x;\Delta) & \triangleq \inf _{\mathbf{v}\in \mathcal{B}^{s,\Delta}} \sup _{\mathbf{u}\in \mathcal{A}^{s,\Delta}} J(s,x;\mathbf{u}, \mathbf{v}) \leq 
V_M^+(s,x;\Delta)\triangleq \inf _{\mathbf{v}\in \mathcal{B}_M^{s,\Delta}} \sup _{\mathbf{u}\in \mathcal{A}^{s,\Delta}} J(s,x;\mathbf{u}, \mathbf{v}).\end{split} \end{equation}
Again, we could have considered above an additional intermediate layer of value functions, playing Markov strategies vs. feedback strategies.
\begin{Remark} One could prove directly that  the game with discretely restricted strategies (either full, or feedback, or Markov) has a value. In other words, we could directly show that (for the more extreme values) 
$$V_M^-(s,x;\Delta)=V_M ^+(s,x;\Delta).$$
This would be rather long technically, and does not bring much understanding into the continuous limit we care about, but it is certainly doable.  The same goes true for the previous Subsection \ref{deterministic}.
\end{Remark}
The main result of the paper, in the general case is:
\begin{Theorem}\label{main-random} Under all standing assumptions on the state equations, together with 
\begin{enumerate}
\item either $p:[0,T]\rightarrow [0,1]$ continuous or
\item both $p$ and $\sigma$ are $C^2$ with respect to the state variable, and  the derivatives
$p_x, \sigma _x$ and $\sigma _{xx}$  exist and are continuous (in all variables),
\end{enumerate}
the Isaacs-type equation \eqref{eq:Isaacs-p} has a unique bounded continuous viscosity solution $v$, which is  the limit (uniform on compact sets)  of  both $V^-_M(\cdot, \cdot;\Delta)$ and $V^+_M (\cdot, \cdot;\Delta) $ as  $\|\Delta \|\rightarrow 0.$ More precisely, for any $K <\infty$ and $\varepsilon >0$, there exists a $\delta=\delta (K, \varepsilon)$ such that
$$\forall s, \forall |x|\leq   K,\ \ \forall \|\Delta\|\leq \delta$$ we have
$$V_M^+(s,x;\Delta)-\varepsilon \leq v(s,x)\leq V_M^-(s,x;\Delta)+\varepsilon.$$
\end{Theorem}

   \section{Proofs} As mentioned, the proofs are based on a modification of Perron's method introduced in \cite{sirbu-4}. However, the analysis has to be done separately for   Subsections  \ref{deterministic}  and  \ref{random}. The proofs for the random case in Subsection \ref{random} are actually easier,  so we will present them more succinctly.
 The notations and definitions  needed for the proofs will refer to the corresponding sub-section. More precisely, we are going to call semi-solutions two different objects in Subsection \ref{deterministic} and in Subsection \ref{random}.

\subsection{Perron  method for deterministic rule of priority}\label{proofs-deterministic}
We start by defining asymptotic stochastic sub and super-solutions for the problem at hand. Before we do that, we need to concatenate strategies (at least the simplest of them, Markov strategies).
\begin{Definition} Fix $s$.  Let $\Delta:t_0<\dots <t_n  \in \mathcal{D}(s)$ and $\xi=(\xi _1, \dots , \xi_n)  \in \mathcal{M}(\Delta)$. Let now $\mathbf{v}\in \mathcal{B}^{s,\Delta, \xi}_M$ be a Markov strategy for player $v$. Fix $t_k\in\Delta$.  The time grid
$\Delta ':  t_k<\dots <t_n=T$ belongs to $\mathcal{D}(t_k)$. Consider now a Markov strategy $\mathbf{v}'\in \mathcal{B}^{t_k, \Delta ', \xi '}_M$  where the marks $\xi'$ are the restriction of the sequence of marks $\xi$ over the interval $[t_k,T]$, i.e.
$$\xi '=(\xi_i)_{i=k+1,\dots,n}.$$
We define the concatenation of $v$ and $v'$ at time $t_k$ over $\Delta, \xi$ by using $v$ up to $t_k$ and $v'$ after that. In other words, we define
$${\bf v\otimes_{t_k} v' }=( v\otimes_{t_k} v'_i)_{i=1,\dots,n},$$
where
$$v\otimes_{t_k} v'_i\triangleq \left \{ \begin{array}{ll}v_i, & 1\leq i \leq k\\
v'_i,  &  k+1\leq i\leq n. \end{array}\right.$$

\end{Definition}
Note that this fits the definition of Markov strategies, if the point $t_k$ is considered as a part of the sub-grid $(r_i)_{i=1,\dots, I(n)}$.
The definitions of asymptotic  sub and super-solutions  below need many quantifiers, but the intuition  behind them is rather clear.

\begin{Definition}[Asymptotic super-solutions]\label{def-super}  A function $w:[0,T]\times \mathbb{R}^d\rightarrow \mathbb{R}$ is called an asymptotic super-solution for the Isaacs equation \eqref{eq:Isaacs-p} if
\begin{enumerate}
\item it is continuous and bounded
\item satisfies $w(T,x)\geq g(x)$ $\forall x\in \mathbb{R}^d,$
\item There exists   a gauge function $\varphi _w :(0,\infty)\rightarrow (0,\infty)$, 
$\lim _{\varepsilon \searrow  0}\varphi  _w(\varepsilon)=0,$
and a constant $C_w>0$ (both possibly depending on $w$) such that;
 $$ \forall \  s\in [0,T], \ \Delta \in \mathcal{D}(s),\ \ \ \xi \in \mathcal{M}(\Delta),\ \ \textrm{and}\ \ \forall t_k\in \Delta,$$ there exist  measurable functions
 $b_w:\mathbb{R}^d\rightarrow V,\ \ \beta_w :\mathbb{R}^d \times U\rightarrow V$ with the property that  (with the notation $\Delta ', \xi'$ as in the above lemma), the Markov strategy   (for the player $v$  starting at $t_k$) $\hat{\mathbf{v}}\in \mathcal{B}_M^{t_k, \Delta ',  \xi'}$ 
defined by $$ \hat{v}_{i}= (b_w (X_{t_k}), \beta_w (X_{t_k}, u),\ \ \ i=k+1, \dots, n,$$
has the property that, 
$\forall u\in \mathcal{A}^{s,\Delta, \xi}, \ \ \forall v\in \mathcal{B}_M^{s, \Delta, \xi},$ $\forall k\leq l\leq n$
 we have 
 \begin{equation}
 \label{eq:supsol}
 \begin{split}
 w(t_k, X^{s, x, \mathbf{u}, \mathbf{v}}_{t_k}) \geq&  \mathbb{E}[w(t_l, X^{s, x, \mathbf{u}, \mathbf{v}\otimes _{t_k} \hat{\mathbf{v}}}_{t_l})| \mathcal{F}_{t_k}] \\- & (t_l-t_k) \left (\varphi  _w(t_l-t_k)+
 C _w \left |\frac { \sum _{i=k+1} ^l (t_i -t_{i-1}) \xi_i }{t_l-t_k}-p(t_k)\right | \right ) .
 \end{split}
 \end{equation}
\end{enumerate}
Denote by $\mathcal{U}$ the class of asymptotic super-solutions.
Comparing   this to the definition of Markov strategies, we see that $(b_w, \beta _w)$  are chosen at $t_k$ and not changed until $t_n$. Only the marks change after $t_k$.\end{Definition}
  \begin{Lemma}\label{supersol-val}Given  $w\in \mathcal{U}$, if $(\Delta ^{n,s}, \xi ^{n,s})$ has uniform asymptotic density $p$ (see the statement of Theorem \ref{main-deterministic}), then 
  $$w(s,x)\geq \limsup _n  V_M^+(s,x;\Delta ^{n,s}, \xi ^{n,s})$$ 
  uniformly in $(s,x)$. 
  More precisely, for any $\varepsilon >0$ there exists $n_0(\varepsilon)$ such that
   \begin{equation}\label{super-sol-uniform}V_M^+(s,x;\Delta ^{n,s}, \xi ^{n,s})-\varepsilon \leq w(s,x),\ \ \forall \ n\geq n_0(\varepsilon),\ \forall s,x.
   \end{equation}
  \end{Lemma}
Proof:   Fix the initial conditions $s,x$ and $\varepsilon>0$. For $n\geq n_0(\varepsilon)$ we have that, inside the partition $$\Delta ^{n,s}=s=t_0<t_1....<t_{N(n)}=T,$$  we can find a sub-division
$$s=r_0<r_1<...r_{I(n)}=T,\ \  \ \ r_i=t_{l(i)},\ \ l(0)=0 <l(1)<...<l(I(n)))=N(n), $$
such that
$$\varphi(r_i-r_{i-1})\leq \varepsilon, i=1,...,I(N), $$
and
$$\left |\frac { \sum _{k=l(i-1)+1} ^{l(i)}(t_k -t_{k-1}) \xi_k }{t_{l(i)}-t_{l(i-1)}}-p(t_{l(i-1)})\right | \leq \varepsilon.$$
Using the definition of super-solution recursively for times   $t_k=t_{l(i-1)}=r_{i-1}$, $i=1, \dots, I(n))$ on the sub-grid, we construct   a $\hat{\mathbf{v}}$ from time $s$ to time $T$ by concatenating the constant strategies $(b^i_w, \beta^i_w)$  (corresponding to time $r_{i-1}$) over the interval $[r_{i-1}, r_i]=[t_{l(i-1)}, t_{l(i)}]$ (remember marks may change over this interval)  such that, in the end, 
$$\hat{v}_l =\left \{ \begin{array}{l r}
\beta^i_w (X|_{[s,r_{i-1}]},u),    & \xi_l=1, \\
b^i_w(X|_{[s, r_{i-1}]}) ,  & \xi_l=0 
\end{array}\right.
\ \ \ l=l(i-1)+1,...,l(i).$$
Note that the strategy $\hat{\mathbf{v}}$ is a Markov strategy.
Implementing this strategy $\hat{\mathbf{v}}$ against \emph{any} strategy $\mathbf{u}$ of the player $u$ yields, according to the very definition of super-solution applied, successively,  for $i=1,....I(n)$ to 
the inequalities
 \begin{equation*}\begin{split}
 w(r _{i-1},X^{s, x, \mathbf{u}, \hat{\mathbf{v}}}_{r_{i-1}})\geq & \mathbb{E}[w(r_i, X^{s, x, \mathbf{u}, \hat{\mathbf{v}}}_{r_i})| \mathcal{F}_{r_{i-1}}]-\\ - &(r_i-r_{i-1}) \left ( \underbrace{\varphi (r_i-r_{i-1})}_{\leq \varepsilon}+
 C \underbrace{\left |\frac { \sum _{k=l(i-1)+1} ^{l(i)}(t_k -t_{k-1}) \xi_k }{r_i-r_{i-1}}-p(r_{i-1})\right | }_{\leq \varepsilon}\right ) .
 \end{split}
 \end{equation*}
 Since $w(T,x)\geq g(x)$, we conclude that
 $$w(s,x)\geq \mathbb{E}[g(X^{s, x, \mathbf{u}, \hat{\mathbf{v}}}_T)]-(T-s)(1+C)\varepsilon, $$ 
 for all $\mathbf{u}$, if $n\geq n_0(\varepsilon)$, finishing the proof of the lemma. $\Diamond$
 
 The proof of the next lemma is not obvious
 but its proof is very similar to the corresponding result in \cite{sirbu-4}.
 \begin{Lemma} \label{closed}The set of asymptotic super-solutions is closed under minimum, i.e. 
 $$v,w\in \mathcal{U} \rightarrow v\wedge w\in \mathcal{U}.$$
 \end{Lemma}
 Next proposition is the main technical result of the section. Before it is presented, notice the obvious fact that  $\mathcal{U}\not= \emptyset$, because the cost function $g$ is bounded.
 \begin{Proposition} \label{perron-inf}Define
 $$v^+\triangleq \inf _{w\in \mathcal{U}} w\geq \limsup _n  V_M^+(s,x;\Delta ^n, \xi ^n).$$
 Then $v^+$ is a USC viscosity sub-solution of the Isaacs equation \eqref{eq:Isaacs-p}.
  \end{Proposition}
  
 Proof: To begin with, we use \cite[Proposition 4.1]{bs-1} together with the previous Lemma \ref{closed} above to conclude that there exist a countable sequence $w_n\in \mathcal{U}$ such that $w_n
\searrow v^+$. We need to treat separately the (parabolic) interior viscosity sub-solution property for $v^+$ and the terminal condition $v^+(T,\cdot)\leq g$ (which is, actually, a equality).

\begin{enumerate}\item {\bf Interior sub-solution property:} Consider a  smooth function $\psi$ which touches $v^+$ locally strictly above at some $(t_0, x_0)\in [0,T)\times \mathbb{R}^d$.  Assume, by contradiction, that
$$\psi_t(t_0, x_0)+H ^p (t_0,x_0,\psi_x(t_0,x_0),\psi_{xx}(t_0,x_0))<-\varepsilon \ \ \textrm{for \ some~} \varepsilon >0.$$
With the notation 
$$L^{t,x,u,v}v=b(t,x,u,v)\cdot v_x (t,x)+\frac{1}{2}Tr\left (\sigma(t,x,u,v)\sigma(t,x,u,v)^T v_{xx}(t,x)\right),$$
this means that there exists a $\hat{b}\in V$ as well as $\hat{\beta}:U\rightarrow V$ such that
$$\psi_t(t_0, x_0)+p(t_0)\sup_{u\in U}L^{t_0,x_0,u, \hat{\beta}(u)} \psi + (1-p(t_0) )\sup_{u\in U}L^{t_0,x_0,u, \hat{b}}\psi < -\varepsilon.$$ 
We make the additional notation
$$\hat{p}\triangleq \psi_t(t_0,x_0)+\sup_{u\in U}L^{t_0,x_0,u, \hat{\beta}(u)} \psi ,\ \ \ \hat{q}\triangleq 
\psi_t(t_0)+\sup_{u\in U}L^{t_0,x_0,u, \hat{b}},$$
and  use continuity of $p$ to consider an even smaller $\varepsilon>0$ such that
$$p(t)(\hat{p}+\varepsilon)+(1-p(t))  ( \hat{q} +\varepsilon)\leq -\varepsilon,\ \ \ |t-t_0|\leq 2\varepsilon.$$
We use continuity once again  to find some $\varepsilon '>0$ depending on $\varepsilon$ such that
$$\psi _t(t,x)+ \sup_{u\in U}L^{t,x,u, \hat{\beta}(u)} \psi\leq \hat{p}+\varepsilon, \psi _t(t,x)+\sup_{u\in U}L^{t_0,x_0,u, \hat{b}} \psi\leq \hat{q}+\varepsilon, \ \ |t-t_0|\vee|x-x_0|\leq 2 \varepsilon '.$$
Since  $\psi$ touches locally $v^+$ strictly above, $\psi$ s continuous and $v^+$ is USC, we can find an even smaller $\varepsilon '>0$ and a   $\delta>0$  small enough such that
$$\psi (t,x)-v^+(t,x)\geq 2 \delta  \textrm{~for~} \varepsilon ' \leq |t-t_0|\vee|x-x_0|\leq 2 \varepsilon '.$$
Now, a Dini type argument, identical to the one used in \cite{bs-2} and \cite{bs-3} implies that, for some $n$ large enough we actually have
\begin{equation}
\label{room}\psi (t,x)-w_n(t,x)\geq  \delta  \textrm{~for~} \varepsilon ' \leq |t-t_0|\vee|x-x_0|\leq 2\varepsilon '.
\end{equation}
Now let $\eta <\delta$ and define
$$w (t,x)\triangleq \left \{
\begin{array}{ll}
(\psi (t,x) -\eta)\wedge w_n(t,x), & |t-t_0|\vee|x-x_0|\leq 2\varepsilon '\\
w_n(t,x),& |t-t_0|\vee |x-x_0|>2\varepsilon '.
\end{array}\right.$$
Using \eqref{room} we see that $w(t,x)=w_n(t,x)$ for $|t-t_0|+|x-x_0|\geq \varepsilon '$ and $w$ is continuous everywhere. On top, $w (t_0,x_0)=v^+(t_0,x_0)-\eta$, so we have a contradiction if we can prove that $w\in \mathcal{U}$. In order to do this,  fix $0\leq s\leq T$,  a $\Delta \in \mathcal{D}(s)$ and some $t_k\in \Delta$.  For this  time $t_k$, we  define $b_w:\mathbb{R}^d\rightarrow V$ and $\beta _w :\mathbb{R}^d\times U\rightarrow V$ by 
$$b_w(x)=\left \{ \begin{array}{ll}
b_{w_n} (x) , & w_n(t_k,x)\leq \psi (t_k,x)-\eta,\\
\hat{b}, &  w_n(t_k,x)> \psi (t_k,x)-\eta,\end{array} \right.$$
as well as 
$$\beta_w(x,u)=\left \{ \begin{array}{ll}
\beta_{w_n} (x,u) , & w_n(t_k,x)\leq \psi (t_k,x)-\eta,\\
\hat{\beta}(u), &  w_n(t_k,x)> \psi (t_k,x)-\eta. \end{array} \right.$$
\end{enumerate}
For any fixed sequence of marks $\xi \in \mathcal{M}(\Delta)$,  fix a $\mathbf{u}\in \mathcal{A}^{s,\Delta, \xi}$.  
Denote by 
$$\hat{X}_{\cdot} \triangleq X_{\cdot}^{s,x,\mathbf{u}, \hat{\mathbf{v}}}.$$Now, the proof goes very similarly to the proof in \cite{sirbu-4}.  Recall that we have to  make all estimates at time $t_k$. On the event where the the initial condition   where $w_n$ lies below $\psi-\eta$ at time $t_k$, 
$$A^c=\{w_n(t_k, \hat{X}_{t_k})\leq \psi (t_k, \hat{X}_{t_k})-\eta\} \in \mathcal{F}_{t_k}$$
equation \eqref{eq:supsol} is satisfied with $\varphi =\varphi _{w_n}$ and $C=C_{w_n}$, i.e. for any later $t_l\in \Delta$, $t_k<t_l$ we have 
\begin{equation}\label{part1}
\begin{split}
1_{A^c} & w(t_k, \hat{X}_{t_k})= 1_{A^c} w_n(t_k, \hat{X}_{t_k})  \geq \\ \geq  & 1_{A^c} \mathbb{E}[w_n (t_l, \hat{X}_{t_l})|\mathcal{F}_{t_k}] -1_{A^c}(t_l-t_k) \left (\varphi  _{w_n}(t_l-t_k)+
 C _{w_n} \left |\frac { \sum _{i=k+1} ^l(t_i -t_{i-1}) \xi_i }{t_l-t_k}-p(t_k)\right | \right )\\
  \geq  & 1_{A^c} \mathbb{E}[w (t_l, X^{r, \xi, {\mathbf{u}, \hat{\mathbf{v}}}}_{t_l})|\mathcal{F}_{t_k}] -1_{A^c}(t_l-t_k) \left (\varphi  _{w_n}(t_l-t_k)+
 C _{w_n} \left |\frac { \sum _{i=k+1} ^l(t_i -t_{i-1}) \xi_i }{t_l-t_k}-p(t_k)\right | \right ). \end{split}
\end{equation}
For the other possibilities, i.e.  over the event 
$$A=\{w_n(t_k, \hat{X}_{t_k})>\psi (t_k, \hat{X}_{t_k})-\eta\}$$
we have, with the additional notation 
$$B=\{ |\hat{X}_t-x_0| \leq 2 \varepsilon '\ \ \forall t_k\leq t\leq t_l\},$$
that 
\begin{equation}
1_A \mathbb{P}( B^c|\mathcal{F}_{t_k})\leq 1_A (t_l-t_k)\varphi  (t_l-t_k)\ \ a.s.,
\end{equation}
uniformly (over all $\mathbf{u}$ and all else that matters), for some gauge function $\varphi$. The proof is  similar to \cite{sirbu-4}. On the other hand, over $A$ we can apply It\^o. On this event,  for any subinterval  $[t_{i-1}, t_i]$ ($i=k+1, \dots, l$) we have
2 possible cases:
\begin{enumerate}
\item either $\xi _i=1$, in which case we use the strategy of $u$ against the counterstrategy $\hat{\beta}$
\item or $\xi_i=0$, when we use a counterstrategy for $u$ agains the constant $\hat{b}$.
After It\^o and taking conditional expectation, we get, 
\begin{equation}\label{part2}\begin{split}
1_Aw (t_k, \hat{X}_{t_k})=& 1_A(\psi (t_k, \hat{X}_{t_k})-\eta)= 1_A \mathbb{E}[\psi (t_l, \hat{X}_{t_l})-\eta |\mathcal{F}_{t_k}] -1_A\Sigma  _{i=k+1}^l \mathbb{E}[\int _{t_{i-1}}^{t_i} L (\psi -\eta) d\tau |\mathcal{F}_r] \\
\geq &  1_A \mathbb{E}[1_B w(t_l, \hat{X}_{t_l})|\mathcal{F}_{t_k}] \\+& 1_A  \mathbb{E}[1_{B^c} \psi (t_l, \hat{X}_{t_l})-\eta|\mathcal{F}_{t_k}]-  1_A \Sigma  _{i=k+1}^l \mathbb{E}[\int _{t_{i-1}}^{t_i} L (\psi -\eta) d\tau |\mathcal{F}_{t_k}] \\=& 
1_A \mathbb{E}[w(t_l, \hat{X}_{t_l})|\mathcal{F}_{t_k}] \\
+& 1_A  \mathbb{E}[1_{B^c} \left (\psi (t_l, \hat{X}_{t_l})-\eta- w(t_l, \hat{X}_{t_l}) \right)|\mathcal{F}_{t_k}]- 1_A \Sigma  _{i=k+1}^l \mathbb{E}[\int _{t_{i-1}}^{t_i} L (\psi -\eta) d\tau |\mathcal{F}_{t_k}] .\end{split}\end{equation}
On $A\cap B$ we have estimates on the generator $L$, so
\begin{equation*}
\begin{split}
1_A \Sigma  _{i=k+1}^l \mathbb{E}[1_B \int _{t_{i-1}}^{t_i} L (\psi -\eta) d\tau |\mathcal{F}_{t_k}] 1\leq 1_A \underbrace{\Sigma  _{i=k+1}^l (t_i-t_{i-1}) \{ (\hat{p}+ \varepsilon)1_{\{\xi_i=1\}}+(\hat{q}+\varepsilon)1_{\{\xi_i=0\}}\} }_{S}.
\end{split}
\end{equation*} 
Now, we have to work with the last sum term
\begin{equation*}
\begin{split}S= &(t_l-t_k)\times  \left  \{(\hat{p}+\varepsilon) \Sigma _{i=k+1}^l\frac{t_i-t_{i-1}}{t_l-t_k}1_{\{\xi_i=1\}} 
+(\hat{q}+\varepsilon) \Sigma _{i=k+1}^l\frac{t_i-t_{i-1}}{t_l-t_k}1_{\{\xi_i=0\}}\right \}\\
=&(t_l-t_k)\times \underbrace{ \left  \{(\hat{p}+\varepsilon) p(t_k)
+(\hat{q}+\varepsilon) (1-p(t_k))\right \} }_{\leq 0}\\
+& (t_l-t_k)\times \left  \{(\hat{p}+\varepsilon) \left (\Sigma _{i=k+1}^l\frac{t_i-t_{i-1}}{t_l-t_k}1_{\{\xi_i=1\}} -p(t_k)\right)
 \right. \\
 & \mbox{}\ \ \ \ \ \ \ \ \ \ \ \ \ \ \ + \left. (\hat{q}+\varepsilon) \left (\Sigma _{i=k+1}^l\frac{t_i-t_{i-1}}{t_l-t_k}1_{\{\xi_i=0\}}-q(t_k)\right)\right \}.\end{split}
\end{equation*}
Therefore,
$$S\leq (t_l-t_k)\times (|\hat{p}+\varepsilon| + |\hat{q}+\varepsilon|) \times \left |\Sigma _{i=k+1}^l\frac{t_i-t_{i-1}}{t_l-t_k}1_{\{\xi_i=1\}} -p(t_k)\right |.$$
On the other hand
$$-1_A  \mathbb{E}[1_{B^c} \left (\psi (t_l, \hat{X}_{t_l})-\eta- w(t_l, \hat{X}_{t_l}) \right)|\mathcal{F}_{t_k}] \leq 1_AC(t_l-t_k)\varphi (t_l-t_k),$$ 
and
 $$1_A \Sigma  _{i=k+1}^l \mathbb{E}[1_{B^c}\int _{t_{i-1}}^{t_i} L (\psi -\eta) d\tau |\mathcal{F}_{t_k}]  \leq 1_A  C (t_l-t_k)^2 \varphi (t_l-t_k) \leq 1_A CT (t_l-t_k) \varphi (t_l-t_k),$$
 where $C$ is a uniform bound on $\psi-\eta -w$ and $L(\psi -\eta)$. Such bound exists, if, for example, we make $\psi$ have compact support. This is fine, since we assumed $\psi$ to touch $v^+$ above strictly only in a local sense.
We now finish the proof by choosing
$$\varphi _w=\varphi _{w_n}\vee C(1+T)\varphi,  \ \ C_w=C_{w_n}\vee \Big( |\hat{p}+\varepsilon | + |\hat{q}+\varepsilon| \Big).$$
\item {\bf The terminal condition: $v^+(T,\cdot)\leq g(\cdot)$}
The proof is very similar to the construction for the terminal condition proof in \cite{sirbu-4} but with the exact modifications present in item 1 above. More precisely, once the  terminal condition  is assumed, by contradiction, to be violated, we  use an identical analytic construction  around the exceptional point $x_0$ such that 
$v^+(T,x_0)>g(x_0)$ to the analytic construction  in \cite{sirbu-4}. The remainder of the proof is then based on two very similar estimates to \eqref{part1} and \eqref{part2} as above to finish the proof. 

\end{enumerate}
\begin{Remark} Obviously, we need to formulate precise counterparts    to Definition \ref{def-super}  to define the class of asymptotic sub-solutions $\mathcal{L}$ and prove a counter-part to Lemma \ref{supersol-val}  and Lemma \ref{closed}. In addition, we construct,, as a counterpart to Proposition \ref{perron-inf} for super-solutions the LSC  viscosity  sub-solution of the Isaacs equation \eqref{eq:Isaacs-p}
$$v^- =\sup _{w\in \mathcal{L}}w,$$
such that $v^-\leq v^+.$
\end{Remark}
The proof of Theorem \ref{main-deterministic} is complete once we realise that
\begin{enumerate}
\item $v^-=v^+=v$ the  continuous viscosity solution to the Isaacs equation \eqref{eq:Isaacs-p}. This is the consequence of  a comparison result whose proof is identical to the result in \cite{sirbu} (first a reduction to a  bounded comparison result, then a limit) once we observe that the (time dependent) linear combination
$$H^p (t,x,p,M)=p(t) H^- (t,x,p,M)+(1-p(t))H^+ (t,x,p,M),$$
satisfies (locally in $x$) the structural condition (3.14) in \cite{cil}, needed for comparison. Note that both Hamiltonians $H^-$ and $H^+$ satisfy such condition (this is proved in \cite[page 19]{cil}) and overviewed in \cite{sirbu}. One has to take care of the appropriate sign for $H$ to get exactly condition (3.14) in \cite{cil}.
\item we can extract a sequence $w_n\searrow v\in \mathcal{U}$ and, therefore, by Dini's criterion (since $w_n$ and $v$ are continuous) we have that the convergence is uniform on compacts.
\item we can make an identical argument for sub-solutions  $v_n\nearrow v$ uniform on compacts,
\item we use relation  \eqref{super-sol-uniform} for  a super-solution $w_n$, and its counterpart for a sub-solution $v_n$  such that $w_n$ and $v_n$ are close apart on compacts (from item 2, 3 above, such exist) to obtain the conclusion \eqref{main-det-precise}.
\end{enumerate}
\subsection{Perron  method for random rules of priority}\label{proofs-random}
We basically have to re-define all objects in Subsection \ref{proofs-deterministic} and go over similar proofs. The arguments will, therefore, be a bit repetitive. For this reason, we present them in less detail. However, care must be taken to account for the randomisation of rules of priority. All notation is considered in the context of the model in Subsection \ref{random}.
\begin{Definition}[Asymptotic super-solutions, random rules]\label{def-super-random}  A function $w:[0,T]\times \mathbb{R}^d\rightarrow \mathbb{R}$ is called an asymptotic super-solution for the Isaacs equation \eqref{eq:Isaacs-p} if
\begin{enumerate}
\item it is continuous and bounded
\item satisfies $w(T,x)\geq g(x)$ $\forall x\in \mathbb{R}^d,$
\item there exist a  gauge function $\varphi _w :(0,\infty)\rightarrow (0,\infty)$,  possibly depending on $w$, with 
$\lim _{\varepsilon \searrow  0}\varphi  _w(\varepsilon)=0,$
such that $\forall s\in [0,T]$  and $ \forall s\leq r\leq T$ there exist measurable functions
$$\hat{b}_w:\mathbb{R}^d\rightarrow V,\ \ \hat{\beta}_w :\mathbb{R}^d \times U\rightarrow V$$ such that the   Markov strategy for the $v$ player 
defined by $\hat{\mathbf{v}}=\mathbf{v} (r)=(\hat{b}_w,\hat{\beta}_w) $ played at time $r$ and kept constant until $T$ against any  $\mathbf{u}=(a,\alpha)$  and following any prior strategy $\mathbf{v}$ up to time $r$  has an asymptotic super-martingale property (until the next time on the grid).    

More precisely, $\forall \Delta \in \mathcal{D}(s)$ such that $r=t_k\in \Delta$, for any $\mathbf{v}\in \mathcal{B}^{s,\Delta}_M$ and any $\mathbf{u}\in \mathcal{A}^{s,\Delta}$, if we make the notation
$\hat{X}_{\cdot}\triangleq X^{s,x, \mathbf{u}, \mathbf{v}\otimes _{t_k} \hat{\mathbf{v}}}_{\cdot},$
  then the "local  asymptotic super-martingale property"
  \begin{equation}\label{eq:supsol-rand}(t_k, \hat{X}_{t_k})\geq  \mathbb{E}[w(t_{k+1}, \hat{X}_{t_{k+1}})| \mathcal{F}_{t_k} \vee \sigma (\eta_1,\dots, \eta_{k})]-(t_{k+1}-t_k)\varphi  _w(t_{k+1}-t_k)
 \end{equation}
 holds.
 \end{enumerate}
Denote by $\mathcal{U}$ the class of asymptotic super-solutions.\end{Definition}
In the definition above we implicitly used a rather obvious notation for the concatenation of Markov strategies.
\begin{Remark}\label{rem-coin-toss}
 We recall that  the priority rule  at time $r=t_k$  is decided by the coin toss $\eta _{k+1}$ and actions are kept constant up to $t_{k+1}$. At time $t_k$, after the coin is tossed, the actions of both players are decided and kept constant until time $t_{k+1}$. Therefore, the process $\hat{X}$ satisfies the integral equation
 \begin{equation*}\begin{split} \hat{X}_t= \hat{X}_{t_k}  +   1_{\{\Phi ^{-1}(\eta _{k+1}\leq   p(t_k,\hat{X}_{t_k})\}} & \left [\int _{t_k}^t  b(\tau, \hat {X}_{\tau}, a_{k+1} (\hat{X}_{\cdot}), \hat{\beta} (\hat{X}_{t_k}, a_{k+1} (\hat{X}_{\cdot})))d\tau  \right. \\+&\left.  \int _{t_k}^t  \sigma (\tau, \hat {X}_{\tau}, a_{k+1} (\hat{X}_{\cdot}), \hat{\beta} (\hat{X}_{t_k}, a_{k+1} (\hat{X}_{\cdot})))d W_{\tau} \right] \\
 + 1_{\{\Phi ^{-1}(\eta _{k+1}>   p(t_k,\hat{X}_{t_k})\}} & \left [\int _{t_k}^t  b(\tau, \hat {X}_{\tau}, \alpha_{k+1} (\hat{X}_{\cdot}, \hat{b}(\hat{X}_{t_k})), \hat{b} (\hat{X}_{t_k})))d\tau  \right. \\+&\left.  \int _{t_k}^t  \sigma 
  (\tau, \hat {X}_{\tau}, \alpha_{k+1} (\hat{X}_{\cdot}, \hat{b}(\hat{X}_{t_k})), \hat{b} (\hat{X}_{t_k})))
 d W_{\tau} \right], t_k\leq t\leq t_{k+1},
  \end{split}
 \end{equation*}
Abusing  notation  we say that
 $$\hat{X}_t=1_{\{\Phi ^{-1}(\eta _{k+1}\leq p(t_k,\hat{X}_{t_k})\}}X^{t_k, \hat{X}_{t_k}, a_{k+1}, \hat{\beta}}_t+
 1_{\{\Phi ^{-1}(\eta _{k+1}> p(t_k,\hat{X}_{t_k})\}}X^{t_k, \hat{X}_{t_k}, \alpha _{k+1}, \hat{b}}_t 
 \in \mathcal{F}_t \vee \sigma (\eta_1,\dots, \eta_{k+1}),$$
 for $t_k\leq t \leq t_{k+1}.$  Taking the expectation with respect tot the coin toss $\eta _{k+1}$, the  "local  asymptotic super-martingale property", written in between consecutive times $t_k$ and $t_{k+1}$ becomes 
  \begin{equation} \label{eq:supsol-rand-exp}  \begin{split}
w(t_k, \hat{X}_{t_k})\geq& p(t_k,\hat{X}_{t_k})\mathbb{E}[w(t_{k+1}, X^{t_k, \hat{X}_{t_k}, a_{k+1}, \hat{\beta}}_{t_{k+1}})| \mathcal{F}_{t_k}\vee \sigma (\eta_1,\dots, \eta_{k}) ] \\ +&  
 (1-p(t_k,\hat{X}_{t_k}) )  \mathbb{E}[w(t_{k+1}, X^{t_k,  \hat{X}_{t_k}, \alpha _{k+1}, \hat{b}}_{t_{k+1} }| \mathcal{F}_{t_k}\vee \sigma (\eta_1,\dots, \eta_{k})]-
 (t_{k+1}-t_k)\varphi  _w (t_{k+1}-t_k).
  \end{split}
 \end{equation}
 We prefer to average out the coin toss $\eta_{k+1}$ because, since can write It\^o formula and identify easy the generator.

\end{Remark}
  \begin{Lemma}\label{supersol-val-rand}Let  $w\in \mathcal{U}$. Then  
  $$w(s,x)\geq \limsup _{\|\Delta|\rightarrow 0}  V_M^+(s,x;\Delta )$$ 
  uniformly in $(s,x)$. 
  More precisely, for any $\varepsilon >0$ there exists $\delta _0(\varepsilon)$ such that
   \begin{equation}\label{super-sol-uniform-rand}V_M^+(s,x;\Delta , \xi )-\varepsilon \leq w(s,x),\ \ \forall \ \Delta \in \mathcal{D}(s), \ \ \|\Delta\|\leq \delta_0(\varepsilon),\ \forall s,x.
   \end{equation}
  \end{Lemma}
Proof:   The proof is actually easier for this case, compared to Subsection \ref{proofs-deterministic}, and very much similar to \cite{sirbu-4}. Fix $s,x$. Fix $\varepsilon>0$. For  some $ \delta_0(\varepsilon)$ we have  
$T\times \varphi _w(\delta)\leq \varepsilon$ for $\delta \leq \delta _0 (\varepsilon).$ Fix a partition 
such that $\|\Delta\|\leq \delta _0(\varepsilon)$.
We use the definition of asymptotic super-solution to construct, recursively, over $[t_k, t_{k+1}]$ the   $\hat{b}_{k+1}, \hat{\beta}_{k+1}$  the strategy for player $v$ 
$$\hat{\mathbf{v}}=(\hat{b}_k, \hat{\beta}_k), \ \ k=1,\dots,n.$$
Use this strategy $\hat{\mathbf{v}}$ against \emph{any} strategy $\mathbf{u}$ of the player $u$ and keep the  notation $\hat{X}$ for the resulting state process.  Iterating \eqref{eq:supsol-rand}  for $t_k=t_0, t_1, \dots, t_{n-1}$ and applying it to up to $t_{k+1}=t_1, t_2, \dots, t_n$ we obtain
 \begin{equation*}
 w(t _{k},\hat{X}_{t_{k}})\geq  \mathbb{E}[w(t_{k+1}, \hat{X}_{t_{k+1}})| \mathcal{F}_{t_k} \vee \sigma (\eta_1,\dots, \eta_{k})] - (t_{k+1}-t_{k})  \underbrace{\varphi _w (t_{k+1}-t_{k})}_{\leq \varepsilon/T}.
 \end{equation*}
 Summing the telescoping terms, since $w(T,x)\geq g(x)$, we conclude that
 $w(s,x)\geq \mathbb{E}[g(X^{s, x, \mathbf{u}, \hat{\mathbf{v}}}_T)]-\varepsilon.$  $\Diamond$
 
 The proof of the next lemma is again very similar to a corresponding result in \cite{sirbu-4}.
 \begin{Lemma} \label{closed-rand}The set of asymptotic super-solutions  with random rules is closed under minimum, i.e. 
 $$v,w\in \mathcal{U} \rightarrow v\wedge w\in \mathcal{U}.$$
 \end{Lemma}
Again,  $\mathcal{U}\not= \emptyset$, because the cost  function $g$ is bounded.
 \begin{Proposition} \label{perron-inf-rand}Define
 $$v^+\triangleq \inf _{w\in \mathcal{U}} w\geq \limsup _{\|\Delta\|\rightarrow 0}  V_M^+(s,x;\Delta ).$$
 Then $v^+$ is a USC viscosity sub-solution of the Isaacs equation \eqref{eq:Isaacs-p}.
  \end{Proposition}
  
 Proof:  Use \cite[Proposition 4.1]{bs-1} and  Lemma \ref{closed-rand}  to find a  sequence $w_n\in \mathcal{U}$ such that $w_n
\searrow v^+$. We again treat separately the (parabolic) interior viscosity sub-solution property for $v^+$ and the terminal condition $v^+(T,\cdot)\leq g$ (equality).

\begin{enumerate}\item {\bf Interior sub-solution property:} Consider a  smooth function $\psi$ which touches locally $v^+$ strictly above at some $(t_0, x_0)\in [0,T)\times \mathbb{R}^d$, and assume that
$$\psi_t(t_0, x_0)+H ^p (t_0,x_0,\psi_x(t_0,x_0),\psi_{xx}(t_0,x_0))<-\varepsilon  \ \textrm{for \ some~} \varepsilon >0.$$
With the same  notation  for
$L^{t,x,u,v}$,
this means that there exists a $\hat{b}\in V$ as well as $\hat{\beta}:U\rightarrow V$ such that
$$\psi_t(t_0, x_0)+p(t_0,x_0)\sup_{u\in U}L^{t_0,x_0,u, \hat{\beta}(u)} \psi + (1-p(t_0,x_0) )\sup_{u\in U}L^{t_0,x_0,u, \hat{b}}\psi < -\varepsilon.$$ 
Now we use continuity to find some  smaller $\varepsilon >0$  such that
\begin{equation}\label{eq:estimates} p(r,y) \left (\psi _t(t,x)+\sup_{u\in U}L^{t,x,u, \hat{\beta}(u)} \psi \right) +(1-p(r,y))\left (\psi _t(t',x')+ \sup_{u\in U}L^{t',x',u, \hat{b}} \psi \right)<-\varepsilon, 
\end{equation}
$\forall  |t-t_0|,|x-x_0|, |r-t_0|, |t'-t_0|, |x'-x_0|\leq 2 \varepsilon $. 
Since the function $\psi$ touches locally $v^+$ strictly above, $\psi$ is continuous and $v^+$ is USC, we can find an even smaller $\varepsilon >0$ and a   $\delta>0$  small enough such that
$$\psi (t,x)-v^+(t,x)\geq 2 \delta  \textrm{~for~} \varepsilon  \leq |t-t_0|\vee |x-x_0|\leq 2 \varepsilon .$$
The same Dini type argument,  implies that, for some $n$ large enough we actually have
\begin{equation}
\label{room}\psi (t,x)-w_n(t,x)\geq  \delta  \textrm{~for~} \varepsilon ' \leq |t-t_0| \vee |x-x_0|\leq 2\varepsilon '.
\end{equation}
Fix  $\eta <\delta$ and define
$$w (t,x)\triangleq \left \{
\begin{array}{ll}
(\psi (t,x) -\eta)\wedge w_n(t,x), & |t-t_0|\vee|x-x_0|\leq 2\varepsilon '\\
w_n(t,x),& |t-t_0|\vee |x-x_0|>2\varepsilon '.
\end{array}\right.$$
Using \eqref{room} we see that $w(t,x)=w_n(t,x)$ for $|t-t_0|\vee|x-x_0|\geq \varepsilon '$ and $w$ is continuous everywhere. Since $w (t_0,x_0)=v^+(t_0,x_0)-\eta$ we have a contradiction if we can prove that  $w\in \mathcal{U}$. We do so next.
For the continuous function $w$, for a fixed $s$ and any $s\leq r\leq T$ define $b_w:\mathbb{R}^d\rightarrow V$ and $\beta _w :\mathbb{R}^d\times U\rightarrow V$ by 
$$b_w(x)=\left \{ \begin{array}{ll}
b_{w_n} (x) , & w_n(r,x)\leq \psi (t,x)-\eta,\\
\hat{b}, &  w_n(r,x)> \psi (t,x)-\eta,\end{array} \right.$$
as well as 
$$\beta_w(x,u)=\left \{ \begin{array}{ll}
\beta_{w_n} (x,u) , & w_n(r,x)\leq \psi (t,x)-\eta,\\
\hat{\beta}(u), &  w_n(r,x)> \psi (t,x)-\eta. \end{array} \right.$$
\end{enumerate}
Fix $s\leq r\leq T$.  Fix also some  $\Delta \in \mathcal{D}(s)$ such that $r=t_k\in \Delta$ for some $k$. Fix any $\mathbf{u}\in \mathcal{A}^{s, \Delta}$ and any $\mathbf{v}\in \mathcal{B}_M^{s, \Delta}$.  
From here on we use all the notations that one would use to write the Definition  \eqref{def-super-random} of for the (potential) super-solution $w$.
 If $w_n$ lies below $\psi-\eta$  at time $t_k=r $ i.e. on 
$$A^c=\{w_n(t_k, \hat{X}_{t_k}\leq \psi (  t_k, \hat{X}_{t_k}  )-\eta\}\in \mathcal{F}_r\vee \sigma (\eta_1, \dots, \eta _{k})$$
equation \eqref{eq:supsol-rand} is satisfied with $\varphi =\varphi _{w_n}$ 
\begin{equation}\label{part1-rand}
\begin{split}
1_{A^c} w(t_k, \hat{X}_{t_k})= & 1_{A^c} w_n( t_k, \hat{X}_{t_k}   ) \\ \geq   & 1_{A^c} \mathbb{E}[w_n (t_{k+1}, \hat{X}_{t_{k+1}})|\mathcal{F}_{t_{k}}\vee \sigma (\eta_1, \dots, \eta _{k})] -1_{A^c}(t_{k+1}-t_k) \varphi  _{w_n}(t_{k+1}-t_k)\\
 \geq   & 1_{A^c} \mathbb{E}[w (t_{k+1}, \hat{X}_{t_{k+1}})|\mathcal{F}_{t_{k}}\vee \sigma (\eta_1, \dots, \eta _{k})] -1_{A^c}(t_{k+1}-t_k) \varphi  _{w_n}(t_{k+1}-t_k). \end{split}
\end{equation}
On  the event 
$A=\{w_n( t_k, \hat{X}_{t_k})>\psi ( t_k, \hat{X}_{t_k})-\eta\}$, 
 with the additional notation 
$$B=\{ |\hat{X}_t-x_0| \leq 2 \varepsilon \ \ \forall r=t_k\leq t\leq t_{k+1}\},$$
we can estimate 
\begin{equation}
1_A \mathbb{P}( B^c|\mathcal{F}_{t_k} \vee \sigma (\eta_1, \dots, \eta _{k}) )\leq 1_A (t_{k+1}-t_k)\varphi  (t_{k+1}-t_k)\ \ a.s.,
\end{equation}
uniformly over all $\mathbf{u}$, for some gauge function $\varphi$; cf. \cite{sirbu-4}. On  the set $A$ we can apply It\^o's formula, separately for the two possibilities for the coin toss at time $r=t_k$ explicitly represented in 
Remark \ref{rem-coin-toss}.
\begin{enumerate}
\item either $\Phi ^{-1}(\eta_{k+1})\leq p(t_k, \hat{X}_{t_k})$, in which case we use the strategy of $u=a_{k+1}$ against the counterstrategy $\hat{\beta}$ 
\item or  $\Phi ^{-1}(\eta _{k+1})> t_k, \hat{X}_{t_k})$, when we use a counterstrategy for $u=\alpha_{k+1}$ agains the constant $\hat{b}$.
It\^o and conditional expectation yields, with the abuse of notation
$$a_{k+1}= a_{k+1}(\hat{X}_{\cdot}), \hat{\beta}=\hat{\beta}(  a_{k+1}(\hat{X}_{\cdot})), \alpha_{k+1}=  \alpha_{k+1}(\hat{X}_{\cdot} ,\hat{b}), $$
and the additional simplifying notation
$$\mathcal{G}_k\triangleq \mathcal{F}_{t_k}\vee \sigma (\eta_1, \dots, \eta _{k}), \psi _{\eta}\triangleq \psi- \eta,$$
\begin{equation}\label{part2-rand}\begin{split}
& 1_A w  (t_k, \hat{X}_{t_k})=   1_A\psi _{\eta}(t_k, \hat{X}_{t_k})=\\ &= 1_A 
\mathbb{E}\left [ 1_{\{\Phi ^{-1}(\eta _{k+1})\leq p( t_k, \hat{X}_{t_k})\}} \left (\psi _{\eta} (t_{k+1}, \hat{X}_{t_{k+1}}) -\int _{t_k}^{t_{k+1}} L^{\tau , \hat{X}_{\tau}, a_{k+1}, \hat{\beta}} \psi _{\eta} d\tau \right )|\mathcal{G}_k \right ]  \\
&+  1_A 
\mathbb{E}\left [ 1_{\{\Phi ^{-1}(\eta _{k+1}> p(t_k, \hat{X}_{t_k})\}} \left (\psi _{\eta} (t_{k+1}, \hat{X}_{t_{k+1}})-\int _{t_k}^{t_{k+1}}L^{\tau, X^{r, \hat{X}_\tau }, \alpha _{k+1}, \hat{b} }\psi _{\eta} d\tau \right )| \mathcal{G}_k \right ]\\
&\geq   1_A \mathbb{E}[1_B w(t_k,\hat{X}_{t_k})| \mathcal{G}_k
] + 1_A  \mathbb{E}[1_{B^c} \psi _{\eta} (t_k, \hat{X}_{t_k})| \mathcal{G}_k]-  \\
&-
1_A \underbrace{\mathbb{E} \left [p(t_k,\hat{X}_{t_k}) \int _{t_k}^{t_{k+1}}  L^{\tau , \hat{X}_{\tau}, a_{k+1}, \hat{\beta}} \psi _{\eta} d\tau  +(1-p(t_k,\hat{X}_{t_k}))   \int _{t_k}^{t_{k+1}} L^{\tau, X^{r, \hat{X}_\tau }, \alpha _{k+1}, \hat{b} }\psi _{\eta} d\tau | \mathcal{G}_k
\right]}_{S}.\end{split}\end{equation}
On the event $A\cap B$ we  have estimates \eqref{eq:estimates} on the generator $L$, and also on the probability (conditional) of $B^c$ so 
$$1_AS\leq 0+ C\mathbb{P}[B^c|\mathcal{F}_{t_k} 
\vee \sigma (\eta_1, \dots, \eta _{k})]\leq 1_A (t_{k+1}-t_k)\varphi (t_{k+1}-t_k).
$$On the other hand
$$-1_A  \mathbb{E}[1_{B^c} \left (\psi _{\eta} (t_{k+1}, \hat{X}_{t_{k+1}})-w( (t_{k+1}, \hat{X}_{t_{k+1}})\right)|\mathcal{F}_{t_k}\vee \sigma (\eta_1, \dots, \eta _{k})] \leq 1_AC(t_{k+1}-t_k)-\varphi (t_{k+1}-t_k),$$ 
 where $C$ is a uniform bound on $\psi_{\eta} -w$ and $L \psi _{\eta}$ (such bound exists, if, for example, we make $\psi$ have compact support by multiplication. This is clearly possible since we assume $\psi$ touches $v^+$ above locally strictly.)
We now finish the proof by choosing
$$\varphi _w=\varphi _{w_n}\vee C(1+T)\varphi.$$
\item {\bf The terminal condition $v^+(T,\cdot)\leq g(\cdot)$}.
The proof is exactly as in the counterpart from Subsection \ref{deterministic},  similar to the construction for the terminal condition, but with the exact modifications present in item 1 above, namely similar estimates to \eqref{part1-rand} and \eqref{part2-rand} are obtained.

\end{enumerate}
\begin{Remark}We again formulate precise counterparts    to Definition \ref{def-super-random}  to define the class of asymptotic sub-solutions $\mathcal{L}$ and prove a counter-part to Lemma \ref{supersol-val-rand}  and Lemma \ref{closed-rand}. In addition, we construct,, as a counterpart to Proposition \ref{perron-inf-rand} for super-solutions the LSC  viscosity  sub-solution of the Isaacs equation \eqref{eq:Isaacs-p}
$$v^- =\sup _{w\in \mathcal{L}}w,$$
such that $v^-\leq v^+.$
\end{Remark}
The proof of Theorem \ref{main-random} (in a similar way to Theorem \ref{main-deterministic})  is complete after we show that 
\begin{enumerate}
\item $v^-=v^+=v$, the  continuous viscosity solution to the Isaacs equation \eqref{eq:Isaacs-p}. This is, again,  the consequence of  a comparison result whose proof is identical to the result in \cite{sirbu}. However, one has to pay special attention to checking (locally in space) the structural condition (3.14) in \cite{cil} for the Hamiltonian
$$H^p((t,x,p,M)=p(t,x)H^-(t,x,p,M)+(1-p(t,x))H^+(t,x,p,M).$$
For the case when $p$ depends on time only, this is all the same as in the proof related to the previous Subsection \ref{deterministic}. However, for the case when $p$ depends on $(t,x)$ this is a  bit more complicated.
We remind the reader, for completeness, that the structural condition (3.14) (localized in space) in \cite{cil} reads, as: 
 $$\forall K<\infty,\ \ \ \exists \ \ \omega_K:[0,\infty)\rightarrow [0,\infty)\ \ \  \omega_k(0+)=0,\textrm{~such~that}$$
\begin{equation}
\label{structural}H(t, x, \alpha (x-y), X) -H(t, y, \alpha (x-y), Y)\leq \omega_K(\alpha |x-y|^2+|x-y|),
\end{equation}
for all $|x|, |y|\leq K, \ t \in [0,T]$ and $X,Y$ such that
\begin{equation}\label{condition}
-3 \alpha \left (\begin{array}{cc} I & 0\\
0& I \end{array} \right)
 \leq \left (\begin{array}{cc} X & 0\\0& -Y \end{array} \right) \leq 3 \alpha \left (\begin{array}{cc} I & -I\\
-I& I \end{array} \right).
\end{equation}
To begin with, we see that,  
$$x\rightarrow p(t,x)b(t,x,u,v)$$
is a locally (in state) Lipschitz map, uniformly in $(t,u,v)$. Therefore, the arguments in \cite{cil} page 19 ensure that
$$G(t,x, u,v, p,M)\triangleq p(t,x)b(t,x,u,v)\cdot p,\ G'(t,x, u,v, p,M)\triangleq (1-p(t,x))b(t,x,u,v)\cdot p$$
satisfy \eqref{structural} subject to \eqref{condition} uniformly in $(t,u,v)$ (locally in $x$).
In addition, the map
$$ x\rightarrow p(t,x)\sigma (t,x,u,v)\sigma ^T (t,x,u,v)\geq I\times 0$$ as a matrix-valued function is in $W^{2, \infty } (B_K)$, where $B_K=\{x\in \mathbb{R}^d| |x|\leq K\},$  uniformly in $t,u,v$.
Theorem 5.2.3  page 132 in \cite{MR2190038} yields that 
$$ x\rightarrow \sqrt{p(t,x)\sigma (t,x,u,v)\sigma ^T (t,x,u,v)} \textrm{~and~} x\rightarrow \sqrt{(1-p(t,x))\sigma (t,x,u,v)\sigma ^T (t,x,u,v)}$$
are (matrix valued maps) locally Lipschitz in $x$, uniformly in $(t,u,v)$.  Following the arguments on page 19 in \cite{cil} we obtain that the "Hamiltonians" 
$$E(t,x,u,v,p,M)\triangleq \frac 12 Tr (p(t,x)\sigma (t,x,u,v)\sigma ^T (t,x,u,v)M),$$
$$ E'(t,x,u,v,p,M)\triangleq \frac 12 Tr ((1-p(t,x))\sigma (t,x,u,v)\sigma ^T (t,x,u,v)M),$$
satisfy \eqref{structural} subject to \eqref{condition} (locally in $x$), and all estimates are uniform in $(t,u,v)$. Following  the comments on last part of page 19 in \cite{cil}, our  Hamiltonian
\begin{equation*}
\begin{split}
H^p(t,x,p,M)=& \sup_{u\in U}\inf_{v\in V} \Big \{G(t,x, u,v,p,M)+E(t,x, u,v,p,M)\Big\} \\+&\inf_{v\in V} \sup_{u\in U}\Big \{G'(t,x, u,v,p,M)+E'(t,x, u,v,p,M)\Big\}
\end{split}
\end{equation*}
also satisfies \eqref{structural} subject to \eqref{condition}, for $|x|,|y|\leq K$.
 \item we can extract a sequence $w_n\searrow v\in \mathcal{U}$ and, therefore, by Dini's criterion (since $w_n$ and $v$ are continuous) we have that the convergence is uniform on compacts.
\item we can make an identical argument for sub-solutions  $v_n\nearrow v$ uniform on compacts,
\item we use relation  \eqref{super-sol-uniform} for  a super-solution $w_n$, and its counterpart for a sub-solution $v_n$  such that $w_n$ and $v_n$ are close apart on compacts (from item 2, 3 above, such exist) to obtain the conclusion \eqref{main-det-precise}.
\end{enumerate}
\bibliographystyle{amsalpha} 
\providecommand{\bysame}{\leavevmode\hbox to3em{\hrulefill}\thinspace}
\providecommand{\MR}{\relax\ifhmode\unskip\space\fi MR }
\providecommand{\MRhref}[2]{%
  \href{http://www.ams.org/mathscinet-getitem?mr=#1}{#2}
}
\providecommand{\href}[2]{#2}

\end{document}